\documentclass[11pt,a4paper,oneside]{article}
\usepackage{amssymb,amsmath}

\usepackage[latin1]{inputenc}
\usepackage{times}
\usepackage{amssymb,amsmath}
\usepackage{latexsym}
\usepackage{amsthm}

\newtheorem{theorem}{{\bf Theorem}}
 \newtheorem{proposition}{{\bf Proposition}}
\newtheorem{corollary}{{\bf Corollary}}
\newtheorem{lemma}{{\bf Lemma}}
\newtheorem{remark}{{\bf Remark}}
\newtheorem{example}{{\bf Example}}

\begin{document}

  \title{ Riemannian Geometry on Contact Lie Groups}

\author{Andr\'e Diatta \footnote{ \footnotesize 
 University of Liverpool. Department of Mathematical Sciences. M$\&$O Building, Peach Street, Liverpool, L69 7ZL, UK. adiatta@liv.ac.uk.
\newline\noindent
$\bullet$ The author was partially supported by Enterprise Ireland.
\newline\noindent
$\bullet$ A part of this work was done while the author was
supported by the IST Programme of the European Union
(IST-2001-35443).}}

\date{}
\maketitle

 \begin{abstract}{ \footnotesize
We investigate contact Lie groups having a left invariant Riemannian
or pseudo-Riemannian metric with specific properties such as being
bi-invariant, flat, negatively curved, Einstein, etc. We classify
some of such contact Lie groups and derive some obstruction results
to the existence of left invariant contact structures on Lie groups
 .}
 \end{abstract}

       \section{\bf Introduction-Summary}
A contact Lie group is a Lie group $G$, say of dimension $2n+1$, having a differential $1$-form $\eta^+$ which  is invariant under left translations by the group elements (left invariant, for short) and which satisfies $(d\eta^+)^n\wedge \eta^+ \neq 0$ pointwise over $G$, where $d\eta^+$ is the de Rham differential of $\eta^+$.
If $\mathcal G$ is the Lie algebra of  $G$ and $\eta:=\eta^+_\epsilon$ the value of $\eta^+$ at the unit $\epsilon$ of $G$, then $(\mathcal G,\eta)$ is called a contact Lie algebra.

Such Lie groups may somehow be seen as playing
the odd dimensional version of Lie groups admitting a left invariant
symplectic structure (symplectic Lie groups). Historically, these
latter have been largely studied by many authors, in line with
 numerous interesting problems in geometry and physics, see e.g. \cite{Shapiro69}, \cite{Lichne-Me88},
\cite{Lichnerowicz90}, \cite{Gindikin-Sapiro-Vinberg},
\cite{drinfeld},\cite{Do-Na}, \cite{dardie-Medi1},
\cite{dardie-Medi2}, \cite{BonYaoChu}. But although contact topology
and geometry are  increasingly acquiring a popular interest due to
their numerous applications,
 contact Lie groups have not been so widely explored (see e.g. \cite{Diatta-Contact}).

 Amongst other results in \cite{Diatta-Contact}, we
 gave a method to construct contact Lie algebras of dimension $2n+1$
  with a  trivial centre, unlike the ordinary contactization (central extension) which only produces
   contact Lie algebras with center of dimension $1$. We discussed some applications of such a construction.
  It turns out that the Lie algebra $\mathcal G_1$ of any exact symplectic Lie group $(G_1,d \alpha^+)$
   can be embedded as a codimension $1$ subalgebra of  many non-isomorphic Lie algebras all of which
having a family of contact forms whose restriction (pullback) to
$\mathcal G_1$ coincide with the $1$-form $\alpha.$
 We also gave a full classification of contact Lie algebras of dimension $5$ and exhibited an infinite family of non-isomorphic solvable contact Lie algebras of dimension $7$.

Here we consider contact Lie groups $G$ which display an additional
structure, namely a left invariant Riemannian or pseudo-Riemannian
metric with specific properties such as being bi-invariant, flat,
negatively curved, Einstein, etc.  The reason for considering these
additional structures lies, on the one hand, in the interest in
(pseudo-) Riemannian geometry.  On the other hand, this can also be
motivated by the fact that the relationship between the contact and
the algebraic structures of Lie groups does not, a priori,  show to
be strong enough to ensure certain general consequences or to affect
certain  invariants of Lie groups.
 So one has to consider specific families of Lie groups.
 Indeed, the main obstructions to the existence of left invariant contact structures on odd dimensional
  Lie groups so far known to the author,
  are the non-degeneracy of their Killing form (Theorem  5 of  \cite{Boothby-Wang58})
   and the dimension of their centre (it is readily checked that the centre should have dimension
   $\leq 1$.)

In this paper, using some known results from Riemannian
    Geometry, we classify contact Lie groups (via their Lie algebras)
    having some given Riemannian or pseudo-Riemannian structure, give
    properties  and derive some obstructions to the existence of left invariant contact
     structures on Lie groups. We carry out a comparison with the symplectic
case, whenever we find it interesting.
     Contact Lie groups turn out to exhibit some behaviour different from that
     of symplectic ones: in presence of some `nice' (pseudo-) Riemannian structures,
      contact structures can exist in abundance
     where under the same assumptions in even dimensions symplectic
      ones would be rather rare and vice versa
       (see Remarks \ref{orthogonal-contact-symplectic}, \ref{rmk:directsums},
\ref{rmk:flatsymplectic}, \ref{nonnegative-Ricci-curv}).

  For the present purposes, we only need to use the presence of left invariant contact and some given
   Riemannian or pseudo-Riemannian structures on the same Lie group.
  The actual relationship of such structures between one another as in   \cite{Blair76}, will be
    discussed in a  subsequent paper. Below, we quote some of our
    main results.

    A Riemannian or pseudo-Riemannian structure on a Lie group is said to be bi-invariant
     if it is invariant under both left and right translations, see \cite{Me-Re85},
     \cite{Me85}.
      The Killing forms of semi-simple Lie groups are examples of such bi-invariant structures.
            In Theorem  5 of  \cite{Boothby-Wang58},  W.M. Boothby and H.C. Wang proved, by generalizing
            a result from J.W. Gray  \cite{Gray57}, that the only
            semi-simple Lie groups that carry a left invariant contact structure are those which are
            locally isomorphic to $SL(2)$ or to $SO(3)$. We extend
               such a result to all Lie groups with bi-invariant Riemannian or pseudo-Riemannian structures.

 \medskip
\noindent
   {\bf Theorem} \ref{contact-orthogonal}
{\em Let $G$ be a                Lie group.
     Suppose                (i) $G$ admits a bi-invariant Riemannian or pseudo-Riemannian
      metric and                (ii) G admits a left invariant contact structure.
       Then $G$ is locally isomorphic to $SL(2,\mathbb R)$ or to $SU(2)$.}

       \medskip
This contrasts with the symplectic case, see Theorem
\ref{split-contact-orthogonal} and Remarks
\ref{orthogonal-contact-symplectic}, \ref{rmk:directsums}.

   In his main result of \cite{Blair-flat} (see also  \cite{Blair76}), D.E. Blair proved that a flat
    Riemannian
    metric in a contact manifold $M$ of dimension $\geq 5$, cannot be a contact metric structure
    (see Section \ref{notation} for the definition).  We prove that in the case of contact Lie
groups of dimension $\geq 5$, there is no flat left invariant
Riemannian metric at all, even if such  a metric has nothing to do
with the given contact structure.

 \medskip
\noindent
   {\bf Theorem \ref{Flat}.}   {\em Let $G$ be a Lie group of dimension $\geq 5$.
 Suppose $G$ admits a left invariant contact structure. Then, there is no
  flat left invariant Riemannian metric on $G$.
More precisely, a contact Lie group admits a flat left invariant Riemannian
 metric if and only if it is locally isomorphic to the group $E(2):=\Bbb R^2\rtimes O(2)$ of rigid motions
 of the Euclidian 2-space.}

Again, we have a different scenario for symplectic Lie groups, as
discussed in Remark \ref{rmk:flatsymplectic}.

\medskip
 We give a characterisation of  contact Lie groups  which have
a left invariant Riemannian metric of negative sectional curvature
(Theorem \ref{codimension1AbIdeal} and Remark
\ref{negativelycurved}).

 \medskip
\noindent {\bf Theorem \ref{Einstein-contact}.} {\em Suppose
$(\mathcal H,\partial \alpha)$ is an exact symplectic solvable Lie
algebra that carries a standard Einstein metric.
 Let $\mathcal A$ be the orthogonal complement
of the derived ideal $[\mathcal H,\mathcal H]$, with respect to the
Einstein metric. Then for any symmetric derivation $D\in
Der(\mathcal H)-\{0\}$ commuting with $ad_a$, for all $a\in \mathcal
A$, the semi-direct product Lie algebra $\mathcal G:= \mathcal
H\rtimes \mathbb                Re$ where $[e,x]=D(x)$, $\forall
x\in \mathcal H$, is a contact Lie algebra endowed with an Einstein
metric.}

\medskip
Clearly, an Einstein metric as in Theorem \ref{Einstein-contact}
does not have a `nice' geometric behaviour with respect to the
contact structure, in the sense that it is not K-contact, as shown
below.

 \medskip
\noindent
 {\bf Theorem 6.}
{\bf A.} {\em If $dim(G)\geq 5$, there is no left invariant contact
structure in any of the following cases: (a) $G$ has the property
that  every left invariant metric has a sectional curvature of
constant sign, (b) $G$ is a negatively curved 2-step solvable Lie
group, (c) $G$ has a left invariant Riemannian metric   with
negative sectional curvature, such that the Levi-Civita connection
$\nabla$ and the curvature tensor $R$ satisfy $\nabla R=0$.}
\newline
\noindent {\bf B.} {\em There is no left invariant K-contact
structure (in dimension $>3$) whose underlying Riemannian metric has
a Ricci curvature of constant sign. In particular, there are no
K-contact-Einstein, a fortiori no Sasaki-Einstein, left invariant
structures on Lie groups of dimension $\geq 5$.}

We split the proof of Theorem 6 into different cases discussed in
the relevant sections, see Proposition \ref{allconstantcurvature},
Corollary \ref{locallysymmetric},  and
Proposition \ref{K-contact}.

\medskip

Most of the results here are
valid for complex contact Lie groups. But as the proofs are the
same, we only explicitly treat the real case.
 The results within this paper might also be useful for
Riemannian or Sub-Riemannian (and CR) Geometry.

\vskip 2truemm \noindent
 { \bf Acknowledgement:} The author  would like to thank Prof. A. Medina and Dr J. Woolf for very helpful
  discussions and
 comments and Prof. A. Agrachev for motivations about Contact Sub-Riemannian Geometry and Vision
 Models. 

  \section{Some preliminaries and notation}\label{notation}

 Unless otherwise stated, $<,>$ will always denote the duality pairing between a vector
  space, say $E$, and its dual space $E^*$. If  $G$ is a Lie group of dimension $2n+1$ and
  $\epsilon$ its unit,
   $\mathcal G$ will stand for its Lie algebra identified with the tangent space
    $T_\epsilon G$ to $G$ at $\epsilon$. When there is no ambiguity, the Lie
    bracket is simply denoted by $[,]$ or by $[,]_{\mathcal G}$ otherwise.
For $x\in \mathcal G$, we will write $x^+$ for the left invariant vector field on $G$ with
 value $x=x^+_\epsilon$  at $\epsilon$.  A left invariant differential $1$-form $\eta^+$ on $G$
 is a contact form  if the left invariant
$(2n+1)$-form $(d\eta^+)^n\wedge \eta^+$ is a volume form on $G$,
where $d\eta^+$ is the de Rham differential of $\eta^+$. This is
equivalent to  $(\partial\eta)^n\wedge \eta $ being a volume form in
$\mathcal G$, where $\eta:=\eta^+_\epsilon$ and
$\partial\eta(x,y):=-\eta([x,y])$.  We call $(G,\eta^+)$ (resp.
$(\mathcal G,\eta)$) a contact Lie group (resp. algebra).  The Reeb
vector field is the unique vector field $\xi^+$ satisfying
$d\eta(\xi^+,x^+)=0$, $\forall x^+$ and $\eta^+(\xi^+)=1$. From now
on, we will also usually write $\partial\eta^+$ instead of
$d\eta^+$. The radical of $\partial \eta$ is the nullspace $Rad
(\partial \eta)=\{x\in \mathcal G, ~ \text{ such that } ~
\partial \eta(x,y)=0, ~~ \forall y\in \mathcal G\}.$ This definition
easily generalizes to any linear $p$-form on $\mathcal G.$

A contact metric structure on a contact manifold $(M,\nu)$ is given by a Riemannian
 metric $g$ and a field  $\phi$ of endomorphisms of its tangent bundle such that for all vector fields $X,Y$,
\begin{equation}\label{contact-metric}  d \nu (X,Y)= g(X,\phi(Y)) \text{ and }
 g(\phi(X),\phi(Y)) =g(X,Y) -\nu(X)\nu(Y)
 \end{equation} (see e.g.  \cite{Blair76}).
If in addition the Reeb vector field is a Killing vector field (i.e.
generates a group of isometries) with respect to $g$, then
 $(g,\phi,\nu)$ is termed a K-contact structure on $M$.
\begin{lemma}\label{kernel-radical} (Lemma 5.2.0.1 of  \cite{Diatta2000}). If $\eta$ is a contact form
 in a Lie algebra $\mathcal G$, with Reeb vector $\xi$, then its kernel $Ker(\eta)$ is not
 a Lie subalgebra of $\mathcal G$, whereas the radical $Rad(\partial\eta)= \mathbb R\xi$ of
 $\partial\eta$ is a reductive subalgebra of $\mathcal G$.\end{lemma}
A symplectic Lie group $(G,\omega^+)$ is a Lie group $G$ together with a left invariant
 symplectic form $\omega^+$ (See e.g. \cite{BonYaoChu}, \cite{dardie-Medi2},  \cite{dardie-Medi1},
 \cite{Lichnerowicz90},  \cite{Lichne-Me88}.)

\noindent
 If the symplectic form $\omega^+$ is the differential
$\partial \alpha^+=\omega^+$ of a left invariant differential 1-form
$\alpha^+$, then   $(G,\partial \alpha^+)$ (resp. $(\mathcal
G,\partial \alpha)$ ) is an exact symplectic Lie group (resp. Lie
algebra).

\noindent
 Exact symplectic Lie algebras of dimension $\leq 6$ are all well known,
  a list of those in dimension $4$ is quoted e.g. in  \cite{Diatta2000}.

\section{\bf Contact Lie groups with a bi-invariant (pseudo-) Riemannian metric}
Our aim in this section is to extend a  result on semi-simple
contact Lie groups due to Boothby and Wang (Theorem $5$ of
\cite{Boothby-Wang58}) to all Lie groups with a bi-invariant
Riemannian or pseudo-Riemannian metric.
 A pseudo-Riemannian metric is a smooth field of bilinear symmetric non-degenerate
 real-valued forms.

                In Theorem $5$ of  \cite{Boothby-Wang58}, Boothby and Wang showed,
                 by generalising a result from J.W. Gray \cite{Gray57},
                  that the only contact Lie groups that are semi-simple are those locally
 isomorphic to $SL(2,\Bbb R)$ or to $SU(2)$. Actually,
 semi-simple Lie groups, with their Killing form, are a
  small part of the much wider family of Lie groups with a Riemannian
   or pseudo-Riemannian metric which is bi-invariant, i.e.
    invariant under both left and right translations.
    For   a connected Lie group, the above property is equivalent to the existence
     of a symmetric bilinear non-degenerate scalar form
     $b$ in its Lie algebra $\mathcal G$, such
    that the adjoint representation of $\mathcal G$ lies in the Lie algebra $\mathcal O(\mathcal G,b)$
     of infinitesimal isometries of $b$, in other words the following holds
     true
\begin{eqnarray}\label{eq:ortogonal}
b([x,y],z)+b(y,[x,z])=0, ~~~ \forall x,y,z\in\mathcal G.
\end{eqnarray}
Such Lie groups and their Lie algebras are called orthogonal (see
e.g. \cite{Me-Re85}, \cite{Me-Re93}).
 This is, for instance, the
case of   reductive Lie groups and Lie algebras  (e.g. the Lie
algebra of all  linear maps $\psi: E\to E$ of a finite dimensional
vector space $E$ ),      the so-called oscillator groups with their
bi-invariant Lorentzian metrics (see \cite{Me85}), the cotangent
bundle of any Lie group (with its natural Lie group structure) and
in general any element of the large and interesting family of the
so-called Drinfeld doubles or Manin algebras which appear as one of
the key tools for the study of the so-called Poisson-Lie groups and
corresponding quantum analogs, Hamiltonian systems (see V.G.
Drinfeld \cite{drinfeld}),      etc. It is then natural to interest
ourselves in the existence of left invariant contact structures on
such Lie groups. Here is our main result.

As well as classifying contact Lie groups that can bear a
bi-invariant Riemannian or pseudo-Riemannian metric, Theorem \ref{contact-orthogonal}
below also implies that the existence of such metrics is an
obstruction to that of left invariant contact structures in
dimension $>3$.
     \begin{theorem}\label{contact-orthogonal} Let $G$ be a                Lie group.
     Suppose                (i) $G$ admits a bi-invariant Riemannian or pseudo-Riemannian
      metric and                (ii) G admits a left invariant contact structure.
       Then $G$ is locally isomorphic to $SL(2,\mathbb R)$ or to $SU(2)$.
\end{theorem}

\begin{remark}\label{orthogonal-contact-symplectic} Unlike the contact Lie groups, in any dimension
 $2n\ge 4$,
  there are several non-isomorphic symplectic Lie
               groups $G$  which also have bi-invariant Riemannian
               or pseudo-Riemannian metrics. The
                 underlying left invariant  symplectic form is related to the bi-invariant
metric by a nonsingular derivation of the Lie algebra $Lie(G)$,
hence $G$ must be nilpotent.
\end{remark}

As a direct corollary of Theorem  \ref{contact-orthogonal}, we have
\begin{theorem}\label{split-contact-orthogonal}
 Suppose a Lie algebra $\mathcal G$ splits as a direct sum
 $\mathcal G = \mathcal G_1 \oplus \mathcal G_2$ of two ideals
 $\mathcal G_1$ and $\mathcal G_2$, where $\mathcal G_1$ is an
 orthogonal Lie algebra. Then $\mathcal G$ carries a contact form if and only
 if $\mathcal G_1$ is  $so(3)$ or $sl(2)$ and $\mathcal G_2$ is an exact symplectic
 Lie algebra.
 \end{theorem}
\proof The proof of Theorem \ref{split-contact-orthogonal} is
straightforward. A Lie algebra which is a direct sum of two ideals
has a contact form if and only if one of the ideals has a contact
form and the other one has an exact symplectic form. From Remark
\ref{orthogonal-contact-symplectic}, if $\mathcal G_1$ had a
symplectic form, then it would have been nilpotent, hence the
symplectic form could not have possibly been exact, due to the
existence of a nontrivial center. So the only possibility is
$\mathcal G_1$ having a contact form and hence being either $sl(2)$
or $so(3)$, and $\mathcal G_2$ having an exact symplectic form. \qed

\begin{remark}\label{rmk:directsums}
  (a) ~~ Theorem  \ref{split-contact-orthogonal} implies in particular that if a Lie algebra
   $\mathcal G$ is a direct sum of its Levi (semi-simple) subalgebra $\mathcal G_1$
   and its  radical (maximal solvable ideal) $\mathcal G_2$,
   then $\mathcal G$ carries a contact form if and only if its Levi component is
    $3$-dimensional and  its radical is an exact symplectic Lie algebra.
     This is  a simple way to construct many non-solvable contact Lie algebras in any
      dimension $2n+1$, where $n\geq 1$.
\newline
\noindent
      (b)~~   Recall that the situation is different in
       the symplectic case.  A symplectic Lie group whose Lie algebra splits as a direct
        sum of its Levi subalgebra and its  radical, must be solvable as shown in
        Theorem 10 of  \cite{BonYaoChu}.
        \end{remark}

   As we need a local isomorphism for the proof of Theorem  \ref{contact-orthogonal},
    we can work with Lie algebras.

            Our following lemma is central in the proof of Theorem \ref{contact-orthogonal}.
\begin{lemma}\label{decompostion-ortho-contact}                If    an orthogonal Lie algebra
$(\mathcal G,b)$ has a contact form $\eta$, then   $\mathcal G$ equals its derived ideal   $\mathcal
               G=[\mathcal                G,\mathcal G]$.                Furthermore, there exists
                $\bar x \in\mathcal G $ such that as a vector space
         $\mathcal G=ker(ad_{\bar x}) \oplus Im(ad_{\bar x})$
                and $ker(ad_{\bar x})$ is of dimension 1, hence $ker(ad_{\bar x}) =\mathbb R\bar x .$
   \end{lemma}
\noindent
 \proof
                Let $\mathcal G$ be a Lie algebra, and $b$ a (possibly non-definite)
                 scalar product on it. For $x\in \mathcal G$,
                 denote by $\theta(x)$ the element of $\mathcal G^*$  defined
                by $<\theta(x) , y>:=b(x,y) $  for all $y$ in $\mathcal G$,
                where  $<,>$ is the duality pairing  between $\mathcal G$ and
                  $ \mathcal G^*$. Then $(\mathcal G, b )$ is an orthogonal Lie algebra if
                   and only if its adjoint  and co-adjoint representations are isomorphic
                   via the linear map
                     $\theta: \mathcal G \to \mathcal G^*$ (see \cite{Me-Re93}).
 Suppose $\eta$ is a contact form on $\mathcal G$. There exists $\bar x $ in $\mathcal G$
 such that $\theta (\bar x) = \eta$. Let us denote by $(\mathbb R \bar
 x)^\bot$ the $b$-orthogonal of $\mathbb R \bar x$, that is, $(\mathbb R \bar
 x)^\bot:=\{y\in \mathcal G ~ \text{ s.t. } ~ b(\bar x, y)=0 \}.$
 Then from the definition of $\bar x$, we have $(\mathbb R \bar
 x)^\bot= \ker(\eta ).$
  The differential of $\eta$ is
  $$\partial \eta (x,y)= - <\eta, [x,y]> = - b(\bar x, [x,y])
   = - b( [\bar x, x],y]).$$ The last equality above is due to Equation
   (\ref{eq:ortogonal}).
 This implies in particular that the radical (nullspace)
 $Rad (\partial \eta)$ of $\partial \eta $ equals the kernel $\ker (ad_{\bar x})$ of
 $ad_{\bar x}$.  As $\eta$ is a contact form,  the vector space underlying $\mathcal  G$
 splits as $\mathcal G:=  Rad (\partial \eta)\oplus \ker (\eta)$ and
 $dim( Rad (\partial \eta))=1$, that is
  $$Rad (\partial \eta)=\ker (ad_{\bar x}) = \mathbb R \bar x.$$
  It then follows that $\dim (Im (ad_{\bar x})) = \dim \mathcal G - 1$.
  Given that $ad_{\bar x}$ is an infinitesimal isometry of $b$, then $Im (ad_{\bar x})$  is a
  subspace of $(\mathbb R \bar x)^\bot $
  and as $\dim (\mathbb R \bar x)^\bot = \dim\mathcal G -1$, we finally conclude that
  $$Im (ad_{\bar x})= (\mathbb R \bar x)^\bot = \ker (\eta).$$
  We have proved that $\mathcal G = \ker (ad_{\bar x}) \oplus Im (ad_{\bar x})$
   and  $\ker (ad_{\bar x}) = \mathbb R \bar x.$

   On the other hand, as $\ker (\eta) $ is not a Lie subalgebra of
   $\mathcal G$ (see lemma  \ref{kernel-radical}), there exist
    $x, y\in \ker (\eta) $, such that  $[x, y]$ is
    not in $\ker (\eta) $, and has the form $[x, y] = t\bar x + [\bar x , x' ]$
     where $t\in \mathbb R -\{0\}$ and $x'\in \mathcal G$.
     But then $\bar x = {1\over t}( [x, y] - [\bar x , x' ])$ is in the derived ideal
      $[\mathcal G ,\mathcal G ]$ of $\mathcal G $
       and consequently we have $\mathcal G = [\mathcal G ,\mathcal G ]$.
\newline
\qed
\newline
\noindent
 {\bf Proof of Theorem \ref{contact-orthogonal}.}
            Let $\mathcal G= \mathcal S\oplus \mathcal R$ be the Levi decomposition of
             $\mathcal G$, where  $\mathcal S$ is the Levi  (semi-simple) subalgebra and
              $\mathcal R$ is the maximal solvable ideal of $\mathcal G$.
              The inequality  $\dim(\mathcal S)\geq 3$ follows from
               Lemma  \ref{decompostion-ortho-contact},
 as the equality $\mathcal G=[\mathcal G,\mathcal G]$ implies that $\mathcal S$ is
               non-trivial.
              We are now going to show that $\mathcal G$ is semi-simple.
 \begin{lemma}  \cite{Me-Re85}\label{ideal-orth-alg}
              A subspace $\mathcal J$ of an orthogonal Lie algebra $(\mathcal G,b)$, is an
              ideal of $\mathcal G$ if and only if the centraliser
              $Z_{\mathcal G}(\mathcal J)$$: = \{ x\in \mathcal G$, such that
              $[x,y] =0$, $\forall y \in \mathcal J\}$ of $\mathcal J$ in $\mathcal G$,
              contains the $b$-orthogonal $\mathcal J^\bot$ of $\mathcal J$.
 \end{lemma}

 \noindent
            Lemma                 \ref{ideal-orth-alg} ensures that
            $Z_{\mathcal G}(\mathcal R)$ contains  $\mathcal R^\bot$ and hence
$$\dim(Z_{\mathcal G}(\mathcal R))\geq \dim(\mathcal R^\bot)=\dim(\mathcal G)-\dim(\mathcal R)
=\dim (\mathcal S)\geq 3.$$
 If the element $\bar x$ as in Lemma
\ref{decompostion-ortho-contact} was in $\mathcal R$, then
$Z_{\mathcal G}(\Bbb R\bar x) = \ker(ad_{\bar x})$ would contain
$Z_{\mathcal G}(\mathcal R)$ and $ \dim( \ker(ad_{\bar x})) \geq 3$.
This would contradict Lemma  \ref{decompostion-ortho-contact}. Of
course $\mathcal R$ being an ideal, implies $ad_{\bar x} \mathcal
R\subset \mathcal R$. Suppose now the restriction $u$ of $ad_{\bar
x}$ to $\mathcal R$ is not injective. Thus, there exists $y_o\neq 0$
in the intersection of $\mathcal R$ and $\ker( ad_{\bar x})$. As
$\bar x$ is not in $\mathcal R$, there exist at least two linearly
independent elements $\bar x, y_o$ in $\ker(ad_{\bar x})$, which
again contradicts Lemma   \ref{decompostion-ortho-contact}. So $u$
is injective and
 $$\mathcal R=u(\mathcal R) \subset Im(ad_{\bar x}) = (\mathbb R \bar x)^\bot
 .$$
  Now the
inclusion $\mathcal R \subset (\mathbb R \bar x)^\bot$ implies $
\mathbb R \bar x \subset \mathcal R^\bot \subset Z_{\mathcal
G}(\mathcal R) $ which mean that $\bar x$ commutes with every
element of $\mathcal R$  and hence this latter is a subset of
$\ker(ad_{\bar x})$. We conclude that $\mathcal R$ is zero, as it is
contained in both $Im(ad_{\bar x})$ and  $\ker(ad_{\bar x})$. So
$\mathcal G$ is semi-simple. But theorem  5 of \cite{Boothby-Wang58}
asserts that the only semi-simple Lie algebras with a contact
structure are $sl(2,\mathbb R)$ and $so(3)$. \qed

\section{Flat Riemannian metrics in Contact Lie Groups}
In his main result of \cite{Blair-flat} (see also \cite{Blair76}),
Blair proved that a contact manifold of dimension $\geq 5$ does not
admit a flat contact metric, ie a metric satisfying the condition
(\ref{contact-metric})  whose sectional curvature vanishes.
Below, we prove that in the case of contact Lie groups of dimension
$\geq 5$, there is no flat left invariant metric at all, even if
such  a metric has nothing to do with the given contact structure.
\begin{theorem}\label{Flat} Let $G$ be a Lie group of dimension $\geq 5$.
 Suppose $G$ admits a left invariant contact structure. Then, there is no
  flat left invariant Riemannian metric on $G$.
\end{theorem}
The following complete classification of contact Lie groups which carry a flat left invariant
metric is a direct consequence of Theorem
  \ref{Flat}.
                 \begin{corollary}A contact Lie group admits a flat left invariant Riemannian
 metric if and only if it is locally isomorphic to the group $E(2):=\Bbb R^2\rtimes O(2)$ of
 rigid motions of the                Euclidian 2-space.
                  \end{corollary}

     \begin{remark}\label{rmk:flatsymplectic}Unlike contact Lie groups  which cannot display
      flat left invariant metrics in
     dimension $>3  $ (Theorem  \ref{Flat}),
      we have again a different scenario for symplectic Lie groups.
      At each even dimension there are several non-isomorphic symplectic Lie groups with
      some flat left invariant metric (see Theorem 2 of
       Lichnerowicz \cite{Lichnerowicz90},  Theorem 2.2 of  \cite{dardie-Medi2}).
       \end{remark}

\noindent
 {\bf  Proof of Theorem \ref{Flat}.} Let $G$
be a connected Lie group of dimension $m$, with a left invariant
Riemannian metric $g$. We denote again by $g$ the scalar product $g_\epsilon$ induced by $g$
 on the Lie algebra $\mathcal G$ of $G$. Then $g$ is flat if and only if its
Levi-Civita connection $\nabla$ defines a homomorphism
$\rho:x\mapsto\rho(x):=\nabla_x$ from the Lie algebra $\mathcal G$
of $G$ to the Lie algebra $ \mathcal O(m)$                consisting
of all skew-adjoint linear maps from $\mathcal G$ to itself.
               This allows Milnor (Theorem  1.5 of  \cite{Milnor76}) to establish that
             $(G,g)$ is flat if and only if $\mathcal G$ splits as  a $g$-orthogonal
             sum $\mathcal G = A_1\oplus A_2$ of a commutative ideal  $A_1:=ker(\rho)$ and
             a commutative subalgebra $A_2$   acting on $A_1$ by skew-adjoint
               transformations obtained by  restricting each $\rho (a)$ to $A_1$,
               for all $a\in A_2$. Let  $\rho$ stand again for such an action of
                $A_2$ on $A_1$ and $\rho^*$ the corresponding contragradiente action of $A_2$ on
                 the dual space $A_1^*$ of $A_1$ by $\rho^* (a)(\alpha):= -
                     \alpha\circ\rho(a)$, for $a\in A_2$ and $\alpha \in A_1^*$.
Namely for $x\in \mathcal A_1$, $a\in \mathcal A_2$, we have
 $$
\rho(a)x:=\nabla_ax=\nabla_ax-\nabla_xa=[a,x].
$$
  Denote $p_i:=dim(A_i)$ the  dimension of $A_i$.
 From the decomposition $\mathcal G=A_1\oplus A_2$, the dual space $\mathcal G^*$ of $\mathcal G$ can be                viewed as $\mathcal G^*=A_2^o\oplus A_1^o$, where $A_i^o$ consists of all linear forms on $\mathcal G$, whose restriction to $A_i$ is identically zero.    All elements of $A^o_1$ are closed forms on $\mathcal G$.   Suppose $\eta=\alpha + \alpha'$ is a contact form on $\mathcal G$, where $\alpha$ is in                $A_2^o\cong  A_1^*$ and $\alpha'$ in $A_1^o\cong A_2^*$. Then $\partial \eta = \partial \alpha $ is given  for all $x,y$  in $A_1$ and all $a,b$ in $A_2$ by $\partial \eta(x,y)=\partial \eta(a,b)=0$,  and
$$\partial \eta(x,a)= -<\alpha, [x,a]>= <\alpha, \rho(a)x>= -<\rho^*(a)\alpha, x>.$$
 Let $m=2n+1$. If $p$ is the dimension of the orbit of $\alpha$ under the action $\rho^*$,
  we can choose linear 1-forms $\alpha_i\in A_{2}^o $  and $\beta_i \in A_{1}^o $ so that
    $\partial \eta$ simply comes to
               $\partial \eta = \displaystyle \sum^{p}_{i=1}                \alpha_i\wedge \beta_i$.
  Due to the property $(\alpha_i\wedge                \beta_i)^2=0$ for each $i=1,...,p$,
   the $2(p+j)$-form $(\partial \eta)^{p+j}$ is identically zero, if $j\geq 1$.
   But obviously we have $p\leq min(p_1,p_2)$.
    Thus as $p_1+p_2=2n+1$, the non-vanishing condition on $(\partial\eta)^n$ imposes that $n=p$
     and either $p_1=p_2+1 = n+1$ or $p_1=p_2-1=n$.
     Hence the dimension of the abelian subalgebra    $\rho(A_2)$ of $\mathcal O(p_1)$ satisfies
      $dim(\rho(A_2)) \geq p \geq p_1-1$. But the maximal abelian subalgebras of $ \mathcal O(p_1)$
      are conjugate to the                Lie algebra of a maximal torus of the compact Lie group
       $SO(p_1)$ (real special orthogonal group of degree $p_1$). It is well known that the dimension of
        maximal tori in $SO(p_1)$ equals ${p_1\over 2}$ if $p_1$ is even, and ${p_1-1\over 2}$ if $p_1$ is odd.
         This is incompatible with the inequality
         $\dim(\rho(A_2))\geq p_1-1$, unless $p_1=2$ and $p_2=1$, hence  $dim(G)=3$.
\qed
\section{Contact Lie Groups with a Riemannian metric of negative curvature }
This section is devoted to the study of contact Lie groups (resp.
algebras) having a left invariant Riemannian metric of negative
sectional curvature. Theorem \ref{codimension1AbIdeal} below,
characterises the more general case of solvable contact Lie algebras
whose derived ideal has codimension $1$. We apply it to the negative
sectional curvature, the locally symmetric and in the $2$-step
solvable cases.
  \begin{theorem} \label{codimension1AbIdeal} \cite{Diatta-Contact}
          (1) If the derived ideal  $\mathcal N:=[\mathcal G, \mathcal G]$ of a solvable contact Lie algebra
           $\mathcal G$ has codimension $1$ in $\mathcal G$, then the following hold. (a) The center $Z(\mathcal N)$  of  $\mathcal N$ has dimension dim$Z(\mathcal N)\le 2$. If moreover dim$Z(\mathcal N)= 2$, then there exists $e\in \mathcal G$, such that $Z(\mathcal N)$ is not an eigenspace of $ad_e$. (b) There is a linear form $\alpha$ on $\mathcal N$ with $(\partial \alpha)^{n-1}\wedge \alpha\neq 0$,     where $dim(\mathcal G)=2n+1$. \newline
          (2) If a Lie algebra $\mathcal G$ has a codimension $1$ abelian subalgebra,
          then $\mathcal G$ has neither a contact form nor an exact symplectic form if dim($\mathcal G)\geq 4$.
\end{theorem}

A detailed proof of Theorem \ref{codimension1AbIdeal} is supplied in \cite{Diatta-Contact}.

\begin{remark} \label{negativelycurved}
 Theorem \ref{codimension1AbIdeal} also characterises
  contact Lie groups with a left invariant Riemannian metric of negative sectional
   curvature. The Lie algebras $\mathcal G$ of such Lie groups are solvable with a
    codimension $1$ derived ideal $[\mathcal G,\mathcal G]$ and there exists $A\in
\mathcal G$ such that the eigenvalues of the restriction of $ad_A$
to $[\mathcal G,\mathcal G]$ have a positive real part. See Theorem
$3$ of \cite{Heintze74}.
\end{remark}

Example \ref{e.g.:codimension1AbIdeal}, exhibits some of the typical
examples of $\mathcal N$  for Theorem \ref{codimension1AbIdeal} (1).
\begin{example}\label{e.g.:codimension1AbIdeal}\cite{Diatta-Contact}
  (i) From a nilpotent symplectic Lie algebra  ($ \mathcal N_o, \omega_o$),
   perform the central extention
$\mathcal N_1 =  \mathcal N_o \times_{ \omega_o} \mathbb R\xi$ using
$\omega_o$, to get  a nilpotent contact Lie algebra with center $Z(
\mathcal N_1)=\mathbb R\xi$. Let a $1$-dimensional Lie algebra
$\mathbb Re_1$ act on $\mathcal N_1$ by a nilpotent derivation $D_1$
with $D_1(\xi)=0$. Denote by $\mathcal N$ the semi-direct product
\begin{eqnarray*}
\mathcal N= \mathcal N_1\rtimes \mathbb Re_1 ~ \text{ so that if } ~
x\in \mathcal N_1 ~ ~\text{ then }~~[e_1,x]:=D_1(x) ~~and ~~
\mathcal N_1 ~\text{ is an ideal of } ~ \mathcal N.\end{eqnarray*}
 Now  we have
\begin{eqnarray*}
Z(\mathcal N) = \mathbb R\xi\oplus \mathbb R(-\bar x+e_1) ~\text{ if
}~ D_1=ad_{\bar x} ~ \text{ for some } ~ \bar x\in \mathcal N_1
 ~\text{ and }~ Z(\mathcal N) = \mathbb R\xi ~ \text{ otherwise.}
\end{eqnarray*}
(ii) Another example is the direct sum $\mathcal N = \mathcal
N'\oplus  \mathcal N''$ of two nilpotent contact Lie algebras
$\mathcal N'$  and $\mathcal N''$, e.g. two Heisenberg Lie algebras
$ \mathcal H_{2p+1}$ and  $\mathcal H_{2q+1}$, thus $Z(\mathcal N)
\subset [\mathcal N,\mathcal N]$ and dim$Z(\mathcal N)$=2.
\end{example}
 In (i) of Example \ref{e.g.:codimension1AbIdeal}, if $D_1=0$, then  $\mathcal N$ is
  the direct sum $\mathcal N= \mathcal N_1\oplus \mathbb Re_1$ of two ideals, namely
   a nilpotent  contact Lie algebra $\mathcal N_1$
   and the line $\mathbb Re_1$. See Example \ref{contact-with-negativecurvature} where
    $\mathcal N_1$ is the Heisenberg Lie algebra $\mathcal H_{2n+1}$.

\begin{example} \label{contact-with-negativecurvature} Let $\mathbb R$ act on the closed connected subgroup
 $$G_n:=\{\sigma=\begin{pmatrix}
1 &x_1&\cdots & x_n&z & 0\\
0&1&0 \cdots & 0&  y_1 & 0\\
\vdots & & & \vdots & \vdots & \vdots
\\
0&\ldots & & 1&y_n&0\\
0&\ldots & & 0&1&0\\
0& \ldots &&0&0&  e^{2\pi iu}
\end{pmatrix}, x_i,y_i,z,u\in\mathbb R\}$$
\noindent
  of $GL(n+3,\mathbb C)$ by
 \begin{eqnarray*}\rho (t) \sigma= (x_1e^{p_1t},
 \cdots, x_ne^{p_nt},y_1e^{(p-p_1)t}, \cdots, y_ne^{(p-p_n)t},ze^{pt},ue^{qt}),\end{eqnarray*}
 \noindent
where $\sigma$ is written as $(x_1, \cdots, x_n,y_1,\cdots,
y_n,z,u)$ for simplicity and $p$, $p_i$, $q$ are constant real
numbers. Let us consider the Lie group $G$ which is the semi-direct
product
$$G=G_n\rtimes_\rho \mathbb R\cong\mathbb R^{2n+1}\times S^1\times
\mathbb R.$$
 The multiplication in $G$ is as follows.
  For every $(\sigma,t):=(x_1,\ldots,x_n,y_1,\ldots,y_n,z,u, t)$,
   $(\tau,s):=(x_1',\ldots,x_n',y_1',\ldots,y_n',z',u', s)$ in $G$,

\begin{eqnarray}
  (\sigma,t)(\tau,s) =
   (x_1+x_1'e^{p_1t}, \ldots,
  x_n+x_n'e^{p_nt},y_1+y_1'e^{(p-p_1)t},\ldots,y_n+y_n'e^{(p-p_n)t},\nonumber\\
   z+z'e^{pt}+\displaystyle\sum_{i=1}^n x_iy_i'e^{(p-p_i)t },u+u'e^{ qt
   }, s+t).
\end{eqnarray}

 If $q\neq p$, then $G$ is a contact Lie group. For instance, for
 every real numbers $k_1,k_2\in \mathbb R$ with $k_1k_2\neq 0$, then
\begin{eqnarray} \eta^+:= k_1e^{-pt}(dz-\displaystyle \sum_{i=1}^n
x_idy_i) +k_2e^{-qt}du
\end{eqnarray}
is a left invariant contact form on $G.$ In fact adding to $\eta^+$
any linear combination (with constant coefficients!) of the left
invariants $1$-forms $e^{-p_it}dx_i$, ~ $e^{(p_i-p)t}dy_i$, ~ $ dt$,
one gets a new left invariant contact form on $G.$

  If moreover
$q>0,p>p_i>0$, $\forall i=1,\ldots, n$ then $G$ has a left invariant
Riemannian metric of negative sectional curvature.

Indeed, the Lie algebra $\mathcal G$ of $G$ is as in Theorem
\ref{codimension1AbIdeal}, see also Example
\ref{e.g.:codimension1AbIdeal}. More precisely, it is a semi-direct
product $\mathcal G=\mathcal N\rtimes \mathbb Re_{2n+2}$ of a
codimension $1$ nilpotent ideal $\mathcal N$ and a $1$-dimensional
Lie subalgebra $\mathbb Re_{2n+2}$. The ideal $\mathcal N$ has a
basis $e_1, \ldots, e_n,e_{n+1},\ldots, e_{2n},e_0, e_{2n+1}$ such
that the semi-direct product $\mathcal N\rtimes \mathbb Re_{2n+2}$
is given by an action of  $\mathbb Re_{2n+2}$ on $\mathcal N$ by
diagonal matrices of the form $\lambda
diag(p_1,\cdots,p_n,p-p_1,\ldots,p-p_n,p,q)$, $\lambda\in \mathbb R$
and $\mathcal N$ itself is a direct sum $\mathcal N:=\mathcal
H_{2n+1}\oplus \mathbb Re_{2n+1}$ of two ideals, $\mathbb Re_{2n+1}$
and the Heisenberg Lie algebra $\mathcal H_{2n+1}=span
(e_1,\ldots,e_{2n}, e_0)$ with centre $\mathbb R e_0$.

 Recall that $G_1$ is
the nilpotent Lie group used by E. Abbena to model the
Kodaira-Thurston Manifold as a nilmanifold, which is symplectic but
not K\"ahlerian. It might be interesting to work out the behaviour
of the extentions to $G$ of the Abbena metric and its relationships
with the contact stuctures on $G$.
\end{example}

 As a direct consequence of Theorem
\ref{codimension1AbIdeal}, we have the following.
\begin{corollary}\label{locallysymmetric} If $dim(G)\geq 5$, then $G$ is not a contact Lie group,
 in any the following cases.
 \newline
{\bf 1.} $G$ is a negatively curved locally symmetric Lie group, i.e. has a left invariant Riemannian metric
   with negative sectional curvature, such that the Levi-Civita connection $\nabla$ and the
    curvature tensor $R$ satisfy $\nabla R=0$. \newline
{\bf 2.} $G$ is a negatively curved 2-step solvable Lie group.
  \end{corollary}
\noindent {\bf Proof.}  (1).  From Proposition  3 of
\cite{Heintze74}, there exists a vector $e$ in the Lie algebra
$\mathcal G$ of $G$ such that $\mathcal G$ splits as a direct sum
$\mathcal G = \mathbb R e\oplus \mathcal A_1\oplus \mathcal A_2$,
where $ \mathcal     N:=\mathcal A_1\oplus \mathcal A_2$ is a 2-step
nilpotent ideal, with derived ideals $[\mathcal G,\mathcal
G]=\mathcal N$ and $[\mathcal N,\mathcal N]= \mathcal A_2$. It
follows that $\mathcal A_2\subset Z(\mathcal N)$.
 But again from  \cite{Heintze74}, $dim (\mathcal A_2)=0,1, 3$, or $7$.
 If $dim (\mathcal A_2)=1$, then $\mathcal G $ has even dimension. The case
  $dim (\mathcal A_2)=0$ corresponds to $\mathcal     N$ being a codimension $1$ abelian
  ideal, which is ruled out, along with the cases $dim(\mathcal     A_2)\geq 3$,
  by Theorem  \ref{codimension1AbIdeal}. So $\mathcal G$ has no contact form.
  The part  (2)  also follows from Theorem
 \ref{codimension1AbIdeal} and Heintze's main result
\cite{Heintze74}, as the derived ideal of the Lie algebra of $G$
must have codimension $1$ and is abelian. \qed

               \begin{proposition}\label{allconstantcurvature} If a Lie group $G$
               has the property that for every left invariant Riemannian metric, the sectional curvature has
                a constant sign, then                $G$ does not carry any left invariant contact
                (or exact symplectic) structure. Moreover, such a Lie group is unique, up to a local
                 isomorphism,  in any dimension.
              \end{proposition}
As a byproduct, the uniqueness result must have another interest
(independent from Contact Geometry)  in the  framework  of
Riemannian Geometry (compare  \cite{Milnor76}, \cite{Nomizu79}).

\vskip 3truemm
 \noindent {\bf Proof of Proposition
\ref{allconstantcurvature}.} From theorem  2.5 of Milnor
\cite{Milnor76} (see also \cite{Nomizu79}), the Lie bracket $[x,y]$
is always equal to a linear combination of $x$ and $y$, for all
$x,y$ in the Lie algebra $\mathcal G$ of such a Lie group. There
exists a well defined real-valued linear map $l$ on $\mathcal G$
such that  $[x,y] = l(y)x - l(x)y$.

               Now identifying the kernel of $l$ with $\mathbb R^n$ and choosing a vector
                $e_1$ satisfying $l(e_1)=1$, allows us to see
                 that all such Lie algebras are actually isomorphic to the sum
                  $\mathbb R^n\oplus \mathbb Re_1$
                  of a codimension $1$ abelian ideal $\mathbb R^n$ and a
complementary $\mathbb Re_1$, where the restriction of $ad_{e_1}$ to
$\mathbb R^n$ is opposite the identity mapping $- id_{\mathbb R^n}$
and $n+1=dim(\mathcal G)$. So any linear form $\alpha$ on $\mathcal
G$, has differential $\partial \alpha=-\alpha\wedge l$. Hence we
have $\partial \alpha\wedge \alpha=0$ and $(\partial \alpha)^p=0$,
$\forall \alpha \in \mathcal G^*$, $\forall p \geq 2$. \qed

\section{\bf Left invariant Einstein metrics on contact  Lie groups }
             As well known, if a connected Lie group $G$ has a left invariant metric with positive Ricci
               curvature, then it must be compact with finite
              fundamental group, see e.g. theorem  2.2 of  \cite{Milnor76}.
                 Thus, Theorem  \ref{contact-orthogonal} ensures that the only
                  Einstein contact Lie groups with a positive
                      Ricci curvature are those locally isomorphic to $SU(2)$.
                       On the other hand, a contact metric structure in a
                        $(2n+1)$-dimensional manifold, is K-contact if and only if the Ricci
                         curvature on the direction of the Reeb vector
                         field $\xi$ is equal to $2n$ (see  \cite{Blair76}).
                          As a direct consequence of this, we have
                \begin{proposition} \label{K-contact}
                 There is no left invariant K-contact structure on Lie groups of
dimension $>3$
                whose underlying Riemannian  metric has a Ricci curvature of constant sign.
                 In particular, there is no K-contact-Einstein, and a fortiori no
                 Sasaki-Einstein, left invariant structure on a Lie group of dimension $\geq 5$.
             \end{proposition}
\begin{remark}\label{nonnegative-Ricci-curv}
Nevertheless, there are contact Lie groups with a left invariant
Riemannian metric of nonnegative Ricci curvature, this is the case
for any $7-$dimensional connected 
 Lie group with Lie algebra
  $\mathbb R^4\rtimes so(3)$ as in Example \ref{e.g.:nonnegativeRiccicurvature}.
  However, the situation is different in the symplectic case. Indeed, as a Lie group with a left
  invariant Riemannian metric of nonnegative Ricci curvature must be unimodular, then from J. Hano
   (see also  \cite{BonYaoChu}) it  is solvable if it admits a left invariant symplectic structure.
   In this case  the metric must be flat (see Lichnerowicz
   \cite{Lichne-Me88}).\end{remark}

\begin{example} \label{e.g.:nonnegativeRiccicurvature} Example of the Lie algebra of a contact Lie
group with a left invariant Riemannian metric of nonnegative Ricci curvature:
$\mathbb R^4 \rtimes so(3)$.

$[e_1,e_2]=e_3,
 [e_2, e_3]=e_1,
[e_3,e_1]=e_2$,

$ [e_1,e_4]=\frac{1}{2}e_7$,
$[e_1,e_5]= \frac{1}{2}e_6$,
$ [e_1,e_6]= \frac{1}{2} e_5$,
$[e_1,e_7]= \frac{1}{2}e_4$,

$[e_2,e_4]= \frac{1}{2}e_5$,
$ [e_2,e_5]= \frac{1}{2}e_4$,
$[e_2,e_6]= \frac{1}{2}e_7$,
$ [e_2,e_7]= \frac{1}{2}e_6$,

$ [e_3,e_4]= \frac{1}{2} e_6$,
  $[e_3,e_5]=\frac{1}{2}e_7$,
 $[e_3,e_6]= \frac{1}{2}e_4$,
$[e_3,e_7]= \frac{1}{2}e_5$,\\
with at least $4$ independent contact forms
$e^*_4,e^*_5,e^*_6,e^*_7$.
\end{example}

Recall that an Einstein metric on a solvable Lie algebra is standard
if the orthogonal complement                of the derived ideal is
an abelian  subalgebra (see e.g.  \cite{Heber98}).
\begin{theorem}\label{Einstein-contact}
Suppose $(\mathcal H,\partial \alpha)$ is an exact symplectic
solvable Lie algebra that carries a standard Einstein metric.
 Let $\mathcal A$ be the orthogonal complement
of the derived ideal $[\mathcal H,\mathcal H]$, with respect to the
Einstein metric. Then for any symmetric derivation $D\in
Der(\mathcal H)-\{0\}$ commuting with $ad_a$, for all $a\in \mathcal
A$, the semi-direct product Lie algebra $\mathcal G:= \mathcal
H\rtimes \mathbb                Re$ where $[e,x]=D(x)$, $\forall
x\in \mathcal H$, is a contact Lie algebra endowed with an Einstein
metric.
              \end{theorem}
\noindent
 \proof
 From  Theorem 2 of \cite{Diatta-Contact}  if $\mathcal G$ is a semi-direct product
  of $(\mathcal H,\partial \alpha)$ and a derivation $D$ of $\mathcal H$, then $\mathcal G$ carries
   a 1-parameter family of contact structures
   $(\eta_t)_{t\in T}$ satisfying $i^*\eta_t=\alpha$, where $i: \mathcal H\to \mathcal G$ is the natural
   inclusion and $T$ is an open nonempty subset of $\mathbb R$. On the other
   hand, from a result of Heber in  \cite{Heber98}, any semi-direct product of a
   standard Einstein Lie algebra $\mathcal H$ by a  symmetric non-trivial derivation
   commuting with $ad_a$, for all  $a\in \mathcal A$, is again a standard Einstein Lie algebra.
\newline
\qed
\newline\noindent
 Theorem  \ref{Einstein-contact} enables us to construct several examples
of contact Lie groups with a left invariant Einstein metric, by
using in particular the so called j-algebras, which are a particular
family of exact symplectic Lie algebra. They also possess an
Einstein metric and  play a central role in the study of the
homogeneous K\"ahler Manifolds and in particular homogeneous bounded
domains. See e.g.  \cite{Do-Na},  \cite{Gindikin-Sapiro-Vinberg},
\cite{Shapiro69}.

 \end{document}